\documentclass[10pt,a4paper,twoside,review]{article}

\usepackage{amssymb}
\usepackage{amsmath}
\usepackage{amsthm}
\usepackage[margin=2.5cm]{geometry}

\def\R{\mathbb{R}}
\def\C{\mathbb{C}}
\def\H{\mathbb{H}}

\def\g{\mathfrak{g}}

\def\p{\mathfrak{p}}

\def\gl{\mathfrak{gl}}

\def\su{\mathfrak{su}}

\def\sp{\mathfrak{sp}}

\def\Dbdle{\mathcal{D}}

\def\Gbdle{\mathcal{G}}
\def\Hbdle{\mathcal{H}}

\def\Qbdle{\mathbb{Q}}

\def\Vbdle{\mathcal{V}}

\def\:{\lrcorner}
\def\#{\sharp}

\def\dsum{\oplus}

\def\prod{\times}

\def\isom{\cong}

\theoremstyle{plain}

\theoremstyle{definition}

\newtheorem*{TheoremA}{Theorem A}

\def\no{\noindent}

{\begin{list}{}%
         {\setlength{\leftmargin}{#1}}%
         \item[]%
}
{\end{list}}

\begin{document}

\title{On the twistor space of a quaternionic contact manifold}

\author{Jesse Alt}

\maketitle

\begin{abstract}
In this note, we prove that the CR manifold which is induced from the canonical parabolic geometry of a quaternionic contact (qc) manifold via a Fefferman-type construction is equivalent to the CR twistor space of the qc manifold defined by O. Biquard.
\end{abstract}




\section{Introduction}

As introduced by O. Biquard in \cite{Biquard}, a \emph{quaternionic contact (qc) manifold} is given by a $4$-tuple $\mathcal{M} = (M,\Dbdle,[g],\Qbdle)$, where $M$ is a manifold of dimension $4n+3$; $\Dbdle \subset TM$ a distribution of co-rank $3$; $[g]$ a conformal class of positive-definite Carnot-Carath\'eodory metrics defined on $\Dbdle$; $\Qbdle$ a rank $3$ sub-bundle of $\mathrm{End}(\Dbdle)$ (all $C^{\infty}$); and where we assume that $\Qbdle$ admits local bases $\{I_1,I_2,I_3\}$ satisfying the quaternion relations (so $I_a^2 = -Id$, $I_1 I_2 = -I_2 I_1 = I_3$) and $\Dbdle$ is given as the kernel of local $1$-forms $\eta^1,\eta^2,\eta^3$, so that the following compatibility relation holds for all $u,v \in \Dbdle$, $a=1,2,3$ and some $g \in [g]$:
\begin{align}
d\eta^a(u,v) = 2g(I_au,v).  \label{qc regularity}
\end{align}
\no In dimension $7$, i.e. for $n=1$, the following integrability condition, due to D. Duchemin \cite{Duch}, will also be assumed: The local $1$-forms $\eta^a$ above may be chosen so that the restrictions of the $2$-forms $d\eta^a$ to $\Dbdle$ form a local oriented orthonormal basis of $\Lambda^2_+\Dbdle^*$, and local vector fields $\xi_1,\xi_2,\xi_3$ (called the \emph{Reeb vector fields} of the $\eta^a$) exist, which satisfy
\begin{align}
\xi_a \: \eta^b = \delta_a^b \,\, \mathrm{and} \,\, (\xi_a \: d\eta^b)_{\vert \Dbdle} = -(\xi_b \: d\eta^a)_{\vert \Dbdle}, \label{qc integrability}
\end{align}
\no for $a,b = 1,2,3$ (in higher dimensions, we always have existence of the Reeb vector fields).\\

A qc structure is naturally defined on the boundary of the rank one symmetric space $\H\mathbf{H}^{n+1} = Sp(1,n+1)/Sp(1)Sp(n+1)$ (the boundary is diffeomorphic to $S^{4n+3}$), and more generally qc structures can be thought of as the natural geometric structures at ``conformal infinity'' of asymptotically symmetric quaternionic-K\"ahler manifolds. Indeed, one of the central results of Biquard's foundational study \cite{Biquard} (Theorem D) says that any real analytic qc manifold $\mathcal{M}$ can be realised as the conformal infinity of a unique asymptotically symmetric quaternionic-K\"ahler metric which is real analytic up to the boundary and defined in a neighbourhood of $\mathcal{M}$.\\

Quaternionic contact structures are the quaternionic analog of Cauchy-Riemann (CR) structures, and there are interesting relations between the two types of geometric structure. An important step in the proof of Biquard's Theorem D is the construction of a natural CR structure on the total space $\mathcal{Z}$ of a $2$-sphere bundle naturally associated to a qc structure $\mathcal{M}$. The space $\mathcal{Z}$ together with this natural CR structure is called the \emph{twistor space} of the qc structure $\mathcal{M}$. (For the construction and proofs of naturality and integrability, cf. II.5 of \cite{Biquard}; we briefly recall the definition in Section 3.)\\

An alternative approach to qc structures is via parabolic geometry: Any qc manifold $\mathcal{M}$ can be canonically identified with a Cartan geometry $(\pi: \Gbdle \rightarrow M,\omega)$ of parabolic type $(G,P)$, where $G \isom Sp(1,n+1)/\{\pm Id\}$ and $P \subset G$ is the parabolic subgroup which is the image under the quotient of the stabiliser in $Sp(1,n+1)$ of a light-like quaternionic line in $\H^{1,n+1}$. That is, $\pi: \Gbdle \rightarrow M$ is a $P$-principal bundle, and $\omega \in \Omega^1(\Gbdle;\g)$ is a Cartan connection of type $(G,P)$. (This is an application of Theorem 3.1.14 of \cite{CSbook}, to which the reader is also referred for background on parabolic geometry; some details of the parabolic structure of a qc manifold are given in Section 2 of \cite{qcweyl}.)\\

Using the Cartan geometry, there is an elegant way to associate a natural CR structure to the qc structure $\mathcal{M}$. Namely, with respect to the inclusion $G \hookrightarrow \tilde{G} := SU(2,2n+2)/\{\pm Id\}$, and for the parabolic subgroup $\tilde{P} \subset \tilde{G}$ which is the quotient of the stabiliser in $SU(2,2n+2)$ of a light-like complex line in $\C^{2,2n+2} \isom \H^{1,n+1}$, we have $\tilde{P} \cap G \subset P$ and $G/(\tilde{P} \cap G) = \tilde{G}/\tilde{P}$. These conditions allow one to execute a Fefferman-type construction (cf. 4.5 of \cite{CSbook} for the general procedure, which includes the application to this specific case in 4.5.5): From $(\pi:\Gbdle \rightarrow M,\omega)$, this construction yields a canonical Cartan geometry $(\tilde{\pi}: \tilde{\Gbdle} \rightarrow \tilde{M},\tilde{\omega})$ of type $(\tilde{G},\tilde{P})$. A Cartan geometry of the latter type (which is also parabolic) is known to induce a partially integrable CR structure of real signature $(4n+2,2)$ on the base space $\tilde{M}$ (some details are recounted in Section 2).\\

Let us refer to the result as the CR Fefferman space of $\mathcal{M}$. In fact, as a by-product of the proof of the main result in \cite{qcfefferman} (cf. Theorem 5.1), the Cartan geometry $(\tilde{\Gbdle},\tilde{\omega})$ of CR type is both normal and torsion-free, and hence (cf. 4.2.4 of \cite{CSbook}) the induced CR structure is integrable. A natural question is as to the relation between this integrable CR structure and the CR twistor space $\mathcal{Z}$ of $\mathcal{M}$. The purpose of this note is to prove that they coincide, confirming the expectation expressed in 4.5.5 of \cite{CSbook}:

\begin{TheoremA} Let $\mathcal{M} = (M,\Dbdle,[g],\Qbdle)$ be a qc manifold (assumed integrable in dimension $7$), and let $(\tilde{\Gbdle} \rightarrow \tilde{M},\tilde{\omega})$ denote the CR Fefferman space induced from the canonical parabolic geometry of $\mathcal{M}$. Then $\tilde{M}$ is naturally identified with the twistor space $\mathcal{Z}$, and the induced CR structures coincide.
\end{TheoremA}

We expect this result to have useful applications for studying the twistor space of a qc manifold, such as computing the Webster scalar curvature for a natural pseudo-hermitian structure on $\mathcal{Z}$ induced by a choice of $g \in [g]$, and these will be developed elsewhere. For now, we mention one immediate corollary of Theorem A: The conformal class of Fefferman metrics of a qc manifold $\mathcal{M}$ (defined on the total space of a natural $S^3$- or $SO(3)$-bundle over $M$, cf. Theorem II.6.1 of \cite{Biquard}) is, up to a finite covering, conformally equivalent to the conformal class of (classical) Fefferman metrics of its twistor space $\mathcal{Z}$ (defined on a natural $S^1$-bundle over a CR manifold, cf. \cite{Fef76}, \cite{Lee}), confirming the expectation expressed in Remark II.6.2 of \cite{Biquard}. This corollary follows because we have proven in \cite{qcfefferman} that the Fefferman-type construction which, using parabolic geometry, induces from $\mathcal{M}$ a natural conformal structure of signature $(4n+3,3)$ is conformally equivalent to the conformal structure defined by Theorem II.6.1 of \cite{Biquard} (this is the equivalence (i) $\Leftrightarrow$ (iv) in Theorem A of \cite{qcfefferman}). On the other hand, carrying out the Fefferman construction of conformal type on $\mathcal{M}$ is obviously equivalent to first carrying out the Fefferman construction of CR type, and then carrying out a Fefferman construction of conformal type on the resulting CR structure. But in \cite{CG06a} it was shown that the result of the latter construction is conformally equivalent, up to a finite covering, to the classical Fefferman metric of a CR manifold.

\section{Background on the flag structures of qc and CR manifolds}

In general, for $G$ a semi-simple Lie group, a parabolic subgroup $P \subset G$ determines an associated $\vert k \vert$-grading of the Lie algebra $\g$ for some $k \in \mathbb{N}$: $\g = \g_{-k} \dsum \ldots \dsum \g_k$ as a vector space, $[\g_i,\g_j] \subset \g_{i+j}$ and $P$ (with Lie algebra $\p = \g_0 \dsum \ldots \dsum \g_k$) consists of the elements in $G$ whose adjoint action preserves the associated filtration $\g = \g^{-k} \supset \ldots \supset \g^{-k}$ (where $\g^i := \g_i \dsum \ldots \dsum \g_k$). The parabolic subgroup has Levi decomposition $P \isom G_0 \ltimes P_+$ where $G_0 \subset P$ is reductive and its adjoint action preserves the grading of $\g$, while $P_+ \subset P$ is a normal, nilpotent subgroup, which is diffeomorphic under the exponential map to $\p_+ := \g^1$, consisting of those elements which strictly increase the grading of elements in $\g$ under the adjoint action. An important object for understanding the underlying geometry on $M$ (called a \emph{flag structure}) which is induced by a Cartan geometry $(\Gbdle \rightarrow M,\omega)$ of parabolic type $(G,P)$, is the bundle $\pi_0: \Gbdle_0 \rightarrow M$, given by $\Gbdle_0 := \Gbdle/P_+$. The filtration of $\g$ induces a filtration of the tangent bundle $TM$ via the isomorphism $TM \isom \Gbdle \times_{\mathrm{Ad}(P)} \g/\p$ (which is general for Cartan geometries), and the Cartan connection $\omega$ identifies the bundle $\Gbdle_0$ as a reduction of the associated graded tangent bundle to $G_0$ (see Chapter 3 of \cite{CSbook}).\\

Now we fix some concepts and notation for the parabolics associated to qc and CR structures, and from here on $(G,P)$ and $(\tilde{G},\tilde{P})$ will denote these fixed parabolic pairs, as indicated in the introduction: First, let $Q$ be the non-degenerate quaternion-hermitian form on $\H^{n+2}$ defined by: $$Q(x) := x_0\overline{x_{n+1}} + \sum_{a=1}^n x_a \overline{x_a} + x_{n+1}\overline{x_0},$$ where we fix the standard ordered basis $\{d_0,\ldots,d_{n+1}\}$ of $\H^{n+2}$ over $\H$ and let $x_i \in \H$ denote the corresponding coordinates of $x$. A calculation yields:
\begin{align}
\g := \sp(Q) &= \{ \left(\begin{array}{ccc}
            a & z & q \\
            \overline{x} & A_0 & -\overline{z}^t \\
            \overline{p} & -x^t & -\overline{a}
                \end{array}\right)  \, \vert \, a \in \H, A_0 \in \sp(n), p,q \in \mathrm{Im}(\H), x,z^t \in \H^n \}, \label{g form}
\end{align}
\no which shows the $\vert 2 \vert$ grading of $\g$ associated to the parabolic subalgebra $$\p := \mathfrak{stab}(\H d_0) = \{ \left(\begin{array}{ccc}
            a & z & q \\
            0 & A_0 & -\overline{z}^t \\
            0 &  0  & -\overline{a}
                \end{array}\right)  \in \g \, \}.$$
We use the form of general elements of $\g$ given by (\ref{g form}) in order to employ a space-saving notation for elements of the specific grading components: E.g., for $p \in \mathrm{Im}(\H)$ we write $[\overline{p}]_{-2} \in \g_{-2}$ to denote the matrix as in (\ref{g form}) with all other entries set to zero; in a similar manner, for $x \in \H^n$ we write $[\overline{x}]_{-1} \in \g_{-1}$ and for $(a,A_0) \in \H \dsum \sp(n) \isom \mathfrak{csp}(1)\sp(n)$ we write $[(a,A_0)]_0 \in \g_0$; etc.\\

Now we let $G := Sp(Q)/\{ \pm Id \}$, which has Lie algebra $\g$, and let $P \subset G$ be the parabolic subgroup (with Lie algebra $\p$) which is the image of the stabiliser in $G$ of $\H d_0$. A further calculation shows that the reductive subgroup preserving the grading components of $\g$ is:
\begin{align}
G_0 = \left\{ \left(\begin{array}{ccc}
            sz & 0 & 0 \\
            0 & A & 0 \\
            0 & 0 & s^{-1}z
                \end{array}\right)  \, \vert \, s \in \R^+, z \in Sp(1), A \in Sp(n) \right\} / \{ \pm Id \}, \label{G 0 form}
\end{align}
\no so $G_0 \isom CSp(1)Sp(n)$ and $P_+ \isom (\H^n)^* \ltimes (\mathrm{Im}(\H))^*$.\\

In Section 2.2 of \cite{qcweyl}, we have given a detailed description of the bundle $\pi_0: \Gbdle_0 := \Gbdle/P_+ \rightarrow M$ in terms of the underlying data $(M,\Dbdle,[g],\Qbdle)$ of a qc manifold, and the explicit action of elements $[(s,z,A)] \in G_0$ on this bundle: A point $u \in \Gbdle_0$ is given by a basis $u = (e_1,\ldots,e_{4n})$ of $\Dbdle_{\pi_0(u)}$ which is symplectic with respect to a metric $g \in [g]$ and a choice of local quaternionic basis $\{I_1,I_2,I_3\}$ of $\Qbdle$ near $\pi_0(u)$. This gives an isomorphism $[u]_{-1}: \Dbdle_{\pi_0(u)} \rightarrow \g_{-1} \isom \H^n$, and hence for $T^{-1}\Gbdle_0 := (T\pi_0)^{-1}(\Dbdle)$ we get a partially-defined $1$-form $\omega_{-1} \in \Gamma(\mathrm{Lin}(T^{-1}\Gbdle_0;\g_{-1}))$ by $\omega_{-1}(\xi) := [u]_{-1}(T_u\pi_0(\xi))$ for $\xi \in T^{-1}_u\Gbdle_0$. By construction, $\omega_{-1}$ is $G_0$-equivariant with respect to the $G_0$-module structure $(\g_{-1},\mathrm{Ad}_{\vert G_0}) \isom (\H^n,\rho_{-1})$, where $\rho_{-1}([(s,z,A)]): \overline{x} \mapsto s^{-1}A(\overline{x})\overline{z}$. In addition, we have a $1$-form $\omega_{-2} \in \Omega^1(\Gbdle_0;\g_{-2})$ which by construction is $G_0$-equivariant with respect to the $G_0$-module structure $(\g_{-2},\mathrm{Ad}_{\vert G_0}) \isom (\mathrm{Im}(\H),\rho_{-2})$, where $\rho_{-2}([(s,z,A)]): \overline{p} \mapsto s^{-2}z \, \overline{p} \, \overline{z}$.\\

Fixing a Carnot-Carath\'eodory metric $g \in [g]$ determines a \emph{scale} for the parabolic geometry $(\Gbdle,\omega)$, and hence a (exact) Weyl structure, i.e. a $G_0$-equivariant section $\sigma: \Gbdle_0 \rightarrow \Gbdle$. Under pull-back via the section $\sigma$, the Cartan connection $\omega$ satisfies: $(\sigma^*\omega_i)_{\vert T^i\Gbdle_0} = \omega_i$, for $i=-1,-2$ as described above. (See \cite{qcweyl}, where the component $\sigma^*\omega_0$ was also computed.) A fixed $g \in [g]$ also determines a complement $\Vbdle \subset TM$ of $\Dbdle$, given as the span of local Reeb vector fields (which is invariant for a fixed $g$).\\

Now let $\tilde{Q}$ be the non-degenerate complex-hermitian form on $\C^{2n+4}$ defined by:
\begin{align*}
\tilde{Q}(y,z) := y_0\overline{y_{n+1}} + \sum_{a=1}^n y_a\overline{y_a} + y_{n+1}\overline{y_0} + z_0\overline{z_{n+1}} + \sum_{a=1}^n z_a\overline{z_a} + z_{n+1}\overline{z_0},
\end{align*}

\no where we identify a vector $y + jz \in \H^{n+2}$ with $(y,z) \in \C^{2n+4}$. We have the standard inclusion $\varphi: \gl(n+2,\H) \hookrightarrow \gl(2n+4,\C)$, given by: $$\varphi: U + jV \mapsto \left(\begin{array}{cc} U & -\overline{V}\\ V & \overline{U}\end{array}\right),$$ and one can verify that this is compatible with the chosen identification $\H^{n+2} \isom \C^{2n+4}$, i.e. that $(U+jV)(y+jz) \simeq \varphi(U+jV)(y,z)$. (In particular, for $\tilde{\p} \subset \tilde{\g} := \su(\tilde{Q})$ the parabolic subalgebra given by $\tilde{\p} := \mathfrak{stab}(\C d_0)$, we have $\varphi^{-1}(\tilde{\p}) \subset \p$.)\\

One can now calculate the decomposition of $\tilde{\g}$ according to the $\vert 2 \vert$-grading associated to $\tilde{\p}$, but to save space we will only give the form of the component $\tilde{\g}_{-1}$, because this is all we need explicitly. We have:

\begin{align*}
\tilde{\g}_{-1} = \left\{ \left(\begin{array}{cc}
                      \left(\begin{array}{ccc} 0 & 0 & 0 \\
                      y & 0 & 0 \\
                      0 & -\overline{y}^t & 0 \end{array}\right) &
                      \left(\begin{array}{ccc} 0 & 0 & 0 \\
                      0 & 0 & 0 \\
                      -\overline{z_+} & -\overline{z}^t & -\overline{z_-} \end{array}\right) \\
                      \left(\begin{array}{ccc} z_- & \,\, 0 \,\, & 0 \\
                      z & 0 & 0 \\
                      z_+ & 0 & 0 \end{array}\right) &
                      0
                \end{array}\right)  \, \vert \,\, y, z \in \C^n, z_-,z_+ \in \C \, \right\}.
\end{align*}
\no (Let us denote an element as above with the row vector $(y,z_-,z,z_+) \in \tilde{\g}_{-1}$.) One verifies that the inclusion $\varphi$ satisfies $\varphi(\p_+) \subset \tilde{\p}$, and that $\varphi(\g^{-1}) \subset \tilde{\g}^{-1}$. Furthermore, if we let $x = x_u + jx_v \in \H^n$, $p = p_u + jp_v \in \mathrm{Im}(\H)$ and $a = a_u + ja_v \in \H$, $A_0 \in \sp(n)$, then we can compute the following formula for the image of elements of $\p$ under the map $\varphi_{-1} = \mathrm{proj}_{\tilde{\g}_{-1}} \circ \varphi$:
\begin{align}
\varphi_{-1}: [\overline{p}]_{-2} &\mapsto (0,0,0,-p_v) \in \tilde{\g}_{-1}; \label{phi -1 projection}\\
\varphi_{-1}: [\overline{x}]_{-1} &\mapsto (\overline{x_u},0,-x_v,0) \in \tilde{\g}_{-1}; \label{phi -1 projection x}\\
\varphi_{-1}: [(a,A&_0)]_0 \mapsto (0,a_v,0,0) \in \tilde{\g}_{-1}. \label{phi -1 projection a}
\end{align}

Letting $\tilde{G}: = SU(\tilde{Q})/\{\pm Id\}$, then the same map gives us an injective homomorphism $\Phi: G \hookrightarrow \tilde{G}$ with differential $\Phi_* = \varphi$. Furthermore, for the (reductive) subgroup $\tilde{G}_0$ of elements which preserve the grading components of the $\vert 2 \vert$-grading of $\tilde{\g}$ associated to $\tilde{\p}$, we have $\tilde{G}_0 \isom (\R_+ \times U(1) \times SU(1,2n+1))/\{ \pm Id \}$.\\

Now let us describe how a Cartan geometry $(\tilde{\pi}: \tilde{\Gbdle} \rightarrow \tilde{M},\tilde{\omega})$ of type $(\tilde{G},\tilde{P})$ induces a (\emph{a priori} partially-integrable) CR structure on the base space: The Cartan connection $\tilde{\omega}$ by definition determines a linear isomorphism $\tilde{\omega}_u: T_u\tilde{\Gbdle} \rightarrow \tilde{\g}$ at each point $u \in \tilde{\Gbdle}$, so in particular this defines a distribution $T^{-1}\tilde{\Gbdle} \subset T\tilde{\Gbdle}$ defined by $T^{-1}_u\tilde{\Gbdle} := \tilde{\omega}_u^{-1}(\tilde{\g}^{-1})$. This defines a distribution $\tilde{\Dbdle} \subset T\tilde{M}$ by letting, for any point $x \in \tilde{M}$, $\tilde{\Dbdle}_x := T_u\tilde{\pi}(T^{-1}_u\tilde{\Gbdle})$ for some $u \in \tilde{\Gbdle}_x$. Since the subspace $\tilde{\g}^{-1} := \tilde{\g}_{-1} \dsum \tilde{\p} \subset \tilde{\g}$ is $\mathrm{Ad}(\tilde{P})$-invariant, and the Cartan connection $\tilde{\omega}$ is $\mathrm{Ad}(\tilde{P})$-equivariant by definition (i.e. $R_p^*\tilde{\omega} = \mathrm{Ad}(p^{-1}) \circ \tilde{\omega}$), it follows that this distribution is well-defined. Also, $\mathrm{rank}_{\R}(\tilde{\Dbdle}) = \mathrm{dim}(\tilde{\g}^{-1}/\tilde{\p}) = \mathrm{dim}(\tilde{\g}_{-1}) = 4n+4$, which shows that $\tilde{\Dbdle}$ is a co-rank $1$ distribution on $\tilde{M}$ (since $\mathrm{dim}(\tilde{M}) = \mathrm{dim}(\tilde{\g}/\tilde{\p}) = 4n+5$).\\

Now let us specify a natural almost complex structure $\tilde{J}$ on $\tilde{\Dbdle}$: Clearly, one can choose a $\tilde{G}_0$-invariant complex structure $J_0$ on $\tilde{\g}_{-1} \isom \C^{2n+2}$ (e.g. scalar multiplication by $-i$), and in fact such a choice is unique up to sign. Since $\tilde{\g}_{-1} \isom \tilde{\g}^{-1}/\tilde{\p}$ as $\tilde{P}$-modules, and $\tilde{P}_+$ acts trivially on $\tilde{\g}^{-1}/\tilde{\p}$, we get a $\tilde{P}$-invariant endomorphism of $\tilde{\g}^{-1}$ from $J_0$ by extending trivially to $\tilde{\p}$, and we'll also denote this by $J_0$. For $x \in \tilde{M}$ and $X \in \tilde{\Dbdle}_x$, choose $u \in \tilde{\Gbdle}_x$ and $\tilde{X} \in T^{-1}_u\tilde{\Gbdle}$ such that $T_u\tilde{\pi}(\tilde{X}) = X$. Then we define $$\tilde{J}(X) := T_u\tilde{\pi}(\tilde{\omega}_u^{-1}(J_0(\tilde{\omega}(\tilde{X})))).$$ Again, equivariance of $\tilde{\omega}$ and $\tilde{P}$-invariance of $J_0$ may be invoked to verify that this definition is proper.\\

\section{Proof of Theorem A}

First let us recall the construction of the twistor space $\mathcal{Z}$ and its CR structure from \cite{Biquard} (we refer also to the exposition in Section 3 of \cite{DavIvMin}): The space $\mathcal{Z} \subset \Qbdle$ is defined fibre-wise, for each point $x \in M$, to be the set of complex structures on $\Dbdle_x$ in $\Qbdle_x$: $$\mathcal{Z}_x := \{ \, I \in \Qbdle_x \, \vert \, I^2 = -Id_{\Dbdle_x} \, \}.$$ This is evidently a $S^2$-bundle over $M$, since any choice of a local quaternionic basis $\{I_1,I_2,I_3\}$ of $\Qbdle$ around $x$ determines an identification of the restriction of $\mathcal{Z}$ to a neighbourhood of $x$ with the endomorphisms $I = a_1I_1 + a_2I_2 + a_3I_3 \in \Qbdle$ such that $a_1^2 + a_2^2 + a_3^2 = 1$.\\

If we fix a choice of Carnot-Carath\'eodory metric $g \in [g]$, then we have a distinguished linear connection $\nabla$ on $M$ (cf. Theorem B, \cite{Biquard}), called the \emph{Biquard connection} of $g$, which induces a horizontal distribution on $\mathcal{Z}$, i.e. we have: $$T_I \mathcal{Z} = Hor^{\nabla}_I (\mathcal{Z}) \dsum Ver_{I}(\mathcal{Z}),$$ where $Ver_{I}(\mathcal{Z}) = T_I(\mathcal{Z}_x)$ is the vertical tangent bundle at $I$ for $I \in \mathcal{Z}_x$. In particular, a choice of $g \in [g]$ determines in this way the horizontal lift of a vector $X \in T_xM$ to $X^{\nabla} \in Hor^{\nabla}_I(\mathcal{Z}) \subset T_I\mathcal{Z}$.\\

A CR distribution $\Hbdle \subset T\mathcal{Z}$ is defined as follows: For $I = a_1I_1 + a_2I_2 + a_2I_3 \in \mathcal{Z}_x$, a corresponding vector $\xi_{I} \in \Vbdle_x \subset T_xM$ is given by letting $\xi_I = a_1 \xi_1 + a_2 \xi_2 + a_3 \xi_3$, where $\xi_1,\xi_2,\xi_3$ are the Reeb vector fields defined locally around $x$ for the unique choice of $1$-forms $\eta^1, \eta^2, \eta^3$ which locally define $\Dbdle$ and are compatible with the local basis $\{I_1, I_2, I_3\}$ and the metric $g$ in the sense of identity (\ref{qc regularity}), and whose existence is assumed in dimension $7$. Declaring $\xi_1, \xi_2, \xi_3$ to be orthonormal, we also have an inner product on $\Vbdle_x$, and a rank $2$ subspace $\xi_I^{\perp} \subset \Vbdle_x$ for any $I \in \mathcal{Z}_x$. Biquard defines: $$\Hbdle_I := (\xi_I^{\perp})^{\nabla} \dsum (\Dbdle_x)^{\nabla} \dsum Ver_I(\mathcal{Z}).$$ Furthermore, an almost complex structure $J^{\mathcal{Z}} \in \Gamma(\mathrm{End}(\Hbdle))$ is defined by letting $J^{\mathcal{Z}}_{\vert (\Dbdle_x)^{\nabla}} = I^{\nabla}$ (the horizontal lift of $I$), and defining the restriction of $J^{\mathcal{Z}}$ to $(\xi_I^{\perp})^{\nabla}$ and to $Ver_I(\mathcal{Z})$ to be the natural complex structures (described explicitly below). Biquard (resp. Duchemin for $n=1$) has proven that the CR structure thus defined is independent of a conformal change of $g \in [g]$, non-degenerate of signature $(4n+2,2)$, and integrable (Theorem II.5.1 of \cite{Biquard}). Once we have identified the twistor space $\mathcal{Z}$ with the CR Fefferman space $(\tilde{\Gbdle},\tilde{\omega})$, these properties follow automatically since $\tilde{\omega}$ is normal and torsion-free.\\

First, let us identify $\tilde{M} \isom \mathcal{Z}$: By definition, $\tilde{M} := \Gbdle/\Phi^{-1}(\tilde{P})$ is the quotient of $\Gbdle$ by the subgroup $\Phi^{-1}(\tilde{P}) \subset P$. Since $P_+ \subset \Phi^{-1}(\tilde{P})$ in our case, and $\Gbdle_0 := \Gbdle/P_+$, we can identify $\tilde{M} = \Gbdle_0/(G_0 \cap \Phi^{-1}(\tilde{P}))$. In fact, in our case the $\R_+$ component of $G_0 \isom CSp(1)Sp(n) = \R_+ \times Sp(1)Sp(n)$ is contained in $\Phi^{-1}(\tilde{P})$, and so for any reduction $\overline{\Gbdle}_0 \hookrightarrow \Gbdle_0$ to the structure group $Sp(1)Sp(n) \subset G_0$, we get an isomorphism $\tilde{M} \isom \overline{\Gbdle}_0/(Sp(1)Sp(n) \cap \Phi^{-1}(\tilde{P}))$. This can be applied, for a fixed choice of $g \in [g]$, to the reduced frame-bundle $\overline{\Gbdle}_0$ consisting, fibre-wise, of those frames $u = (e_1,\ldots,e_{4n})$ of $\Dbdle_x$ which are symplectic with respect to $g$ and some local quaternionic basis $\{I_1,I_2,I_3\}$ of $\Qbdle$.\\

Note that the subgroup $Sp(1)Sp(n) \cap \Phi^{-1}(\tilde{P})$ consists, with respect to the presentation (\ref{G 0 form}), of those elements $[(1,z,A)] \in G_0$ for which $z \in Sp(1)$ is of the form $z = z_0 + z_1i$, i.e. $z \in U(1) \subset Sp(1)$. Thus, $Sp(1)Sp(n) \cap \Phi^{-1}(\tilde{P}) = U(1)Sp(n)$. Also, we see that the subgroup $Sp(1)Sp(n) \cap \Phi^{-1}(\tilde{P})$ equals the stabiliser of the point $i \in S^{2}$ (identifying $S^2$ with the unit imaginary quaternions) under the action $\rho_0: Sp(1)Sp(n) \rightarrow Diff(S^2)$ given by $\rho_0([1,z,A]): q \mapsto z q \overline{z}$.\\

For a point $u \in \overline{\Gbdle}_0$ given by a symplectic basis of $\Dbdle$ with respect to $g$ and $\{I_1,I_2,I_3\}$, we identify $(u,i) \simeq I_1$, $(u,j) \simeq I_2$ and $(u,k) \simeq I_3$. By definition of the $G_0$ action on $\Gbdle_0$ (cf. Section 2.2 of \cite{qcweyl}), this identification is invariant under the right action by $Sp(1)Sp(n)$ on $\overline{\Gbdle}_0 \times S^2$ given by $(u,q).g = (ug,\rho_0(g^{-1})(q))$. In particular, this gives an isomorphism $\overline{\Gbdle}_0 \times_{\rho_0} S^2 \isom \mathcal{Z}$ as $S^2$-bundles over $M$. On the other hand, $\overline{\Gbdle}_0 \times_{\rho_0} S^2 \isom \overline{\Gbdle}_0/(Sp(1)Sp(n) \cap \Phi^{-1}(\tilde{P}))$, since $Sp(1)Sp(n) \cap \Phi^{-1}(\tilde{P})$ is the stabiliser of $i$ under the action $\rho_0$.\\

In summary, this gives us the identification $\tilde{M} \isom \mathcal{Z}$ (and hence a submersion which we'll denote $\wp_0: \overline{\Gbdle}_0 \rightarrow \mathcal{Z}$, given by sending $u \in \overline{\Gbdle}_0$ to the point in $\mathcal{Z}$ identified with $[(u,i)] \in \overline{\Gbdle}_0 \times_{\rho_0} S^2$). From the preceding argument, the following fact about this identification is evident: Fixing a local quaternionic basis $\{I_1,I_2,I_3\}$ of $\Qbdle$ around $x \in M$ and a point $u \in (\overline{\Gbdle}_0)_x$ corresponding to this basis (and to $g \in [g]$, which is fixed throughout), then for any point $I = a_1I_1(x) + a_2I_2(x) + a_3I_3(x) \in \mathcal{Z}_x$ we have $I = \wp_0(ug_I)$ where $g_I = [(1,z_I,Id)]$ for some $z_I \in Sp(1)$ such that $\rho_0(g_I)(i) = z_I i \overline{z_I} = a_1i + a_2j + a_3k$. This will be useful for subsequent calculations.\\

In the next step, we identify the induced CR distribution on $\mathcal{Z}$ corresponding to $\tilde{\Dbdle} \subset T\tilde{M}$. We will abuse notation slightly by writing $\tilde{\Dbdle}_I \subset T_I \mathcal{Z}$ for $I \in \mathcal{Z}$. From the construction of the Fefferman space $(\tilde{\Gbdle},\tilde{\omega})$, we have an inclusion $\iota: \Gbdle \hookrightarrow \tilde{\Gbdle}$ of bundles over $\tilde{M}$, and the Cartan connections are related by $\iota^*\tilde{\omega} = \varphi \circ \omega$. Moreover, the composition $\tilde{\pi} \circ \iota: \Gbdle \rightarrow \tilde{M}$ equals the defining projection $p: \Gbdle \rightarrow \tilde{M} := \Gbdle/\Phi^{-1}(\tilde{P})$. Denoting the induced projection by $\wp: \Gbdle \rightarrow \tilde{M}$, we thus have $\tilde{\Dbdle}_I = T_{\tilde{u}}\wp((\varphi \circ \omega_{\tilde{u}})^{-1}(\tilde{\g}^{-1}))$ for a choice of $\tilde{u} \in \wp^{-1}(I)$. If we denote by $\sigma: \Gbdle_0 \rightarrow \Gbdle$ the $G_0$-equivariant section (Weyl structure) corresponding to the Carnot-Carath\'eodory metric $g \in [g]$, let $u \in \overline{\Gbdle}_0$, $I = a_1I_1 + a_2I_2 + a_2I_3$ and $g_I, z_I$ be as in the preceding paragraph, then $ug_I \in \wp_0^{-1}(I)$ and we have $$\tilde{\Dbdle}_I = T_{ug_I}\wp_0((\varphi \circ \sigma^*\omega)^{-1}(\tilde{\g}^{-1})).$$

As noted in Section 2, we have $\varphi(\g^{-1}) \subset \tilde{\g}^{-1}$. Hence, for $I \in \mathcal{Z}_x$ and any $X \in \Dbdle_x \subset T_xM$, if $\tilde{X} \in T_{ug_i}\overline{\Gbdle}_0$ is any lift of $X$ to the point $ug_I$, then $T_{ug_I}\wp_0(\tilde{X}) \in \tilde{\Dbdle}_I$ (since $\sigma^*\omega(\tilde{X}) \in \g^{-1}$). Also, any vertical (over $M$) tangent vector in $T_{ug_I}\overline{\Gbdle}_0$ projects to $\tilde{\Dbdle}_I$. In particular, $(\Dbdle_x)^{\nabla} \dsum Ver_I(\mathcal{Z}) \subset \tilde{\Dbdle}_I$, since we can take the horizontal lift $X^h$ of any vector $X \in \Dbdle_x$ to $ug_I$ with respect to the Biquard connection form on $\overline{\Gbdle}_0$, which clearly projects to $X^{\nabla} \in T_I\mathcal{Z}$.\\

To show the inclusion $(\xi_I^{\perp})^{\nabla} \subset \tilde{\Dbdle}_I$ (and hence the equality $\Hbdle_I = \tilde{\Dbdle}_I$), we first look closer at the image $\varphi(\g_{-2}) \subset \tilde{\g}$: Namely, one calculates that $\varphi([-j]_{-2}), \varphi([-k]_{-2}) \in \tilde{\g}_{-1} \subset \tilde{\g}^{-1}$. For $I \in \mathcal{Z}_x$ and $u \in \overline{\Gbdle}_0$, $g_I \in Sp(1)Sp(n)$, $z_I \in Sp(1)$ as specified above, we define $J, K \in \mathcal{Z}_x$ by $J := b_1I_1 + b_2I_2 + b_3I_3$ and $K := c_1I_1 + c_2I_2 + c_3I_3$, for $b_1i + b_2j + b_3k := z_I j \overline{z_I}$ and $c_1i + c_2j + c_3k := z_I k \overline{z_I}$. Then $\xi_J, \xi_K \in \Vbdle_x$ span the orthogonal complement of $\xi_I$ in $\Vbdle_x$. By construction, $\omega_{-2}(u)(\tilde{\xi_J}) = [-z_I j \overline{z_I}]_{-2}$ and $\omega_{-2}(u)(\tilde{\xi_K}) = [-z_I k \overline{z_I}]_{-2}$ for $\tilde{\xi_J}, \tilde{\xi_K}$ any lifts of $\xi_J, \xi_K$, respectively, to the point $u$. Using $G_0$-equivariance, we get:
\begin{align*}
\omega_{-2}(ug_I)((R_{g_I})_*(\tilde{\xi_J})) &= (R_{g_I}^*\omega_{-2})(u)(\tilde{\xi_J}) = \mathrm{Ad}(g_I^{-1})(\omega_{-2}(u)(\tilde{\xi_J})) \\
&= \overline{z_I} \, (\omega_{-2}(u)(\tilde{\xi_J})) \, z_I = [-j]_{-2}.
\end{align*}
\no Similarly, $\omega_{-2}(ug_I)((R_{g_I})_*(\tilde{\xi_K})) = [-k]_{-2}$. Thus $(R_{g_I})_*(\tilde{\xi_J}), (R_{g_I})_*(\tilde{\xi_K}) \in (\varphi \circ \sigma^*\omega)^{-1}(\tilde{\g}^{-1}) \subset T_{ug_i}\overline{\Gbdle}_0$ (and hence any lifts of the vectors $\xi_J, \xi_K$ to the point $ug_I \in \overline{\Gbdle}_0$) project via $\wp_0$ to $\tilde{\Dbdle}_I$, so $(\xi_I^{\perp})^{\nabla} \subset \tilde{\Dbdle}_I$.\\

It remains to compute the induced almost complex structure $\tilde{J}$ on $\tilde{\Dbdle}$. For this calculation, we only need to consider the components $(\varphi_{-1} \circ \sigma^*\omega_{\leq})(\tilde{X})$ for $\tilde{X} \in T\overline{\Gbdle}_0$, where $\omega_{\leq} := \omega_{-2} + \omega_{-1} + \omega_0$, since $\varphi(\p_+) \subset \tilde{\p}$. First, note that for any $X \in \Dbdle_x$, the horizontal lifts of $X$ to vectors in $\overline{\Gbdle}_0$ with respect to the Biquard connection and the Weyl connection $\omega_0 = \sigma^*\omega_0$ are the same, so in particular we have $\sigma^*\omega_{\leq}(X^h) = \omega_{-1}(X^h)$ for $X^h$ the horizontal lift via the Biquard connection. This follows from the computation of the Weyl connection with respect to the Weyl structure $\sigma$ induced by $g \in [g]$, cf. Theorem 3.7 of \cite{qcweyl}. On the other hand, let us denote by $\xi_a^h \in T_u\overline{\Gbdle}_0$ the horizontal lift of a Reeb vector field $\xi_a \in \Vbdle_x$ to the point $u$ with respect to the Biquard connection (where $u$ and the basis $\{I_1,I_2,I_3\}$ are related as specified above). Then it follows from the same result that $\sigma^*\omega_{\leq}(\xi_a^h) = \omega_{-2}(\xi_a^h) + \omega_0(\xi_a^h)$ and we have:
\begin{align}
\omega_0(\xi_a^h) = [(\tilde{\mathrm{s}}_g i_a,\omega_A(\xi_a^h))]_0 \in \g_0 \label{omega 0 form}
\end{align}
where $\tilde{\mathrm{s}}_g := \mathrm{scal}/32n(n+2)$ is the rescaled qc scalar curvature of $g \in [g]$ and $i_1 := i, i_2 := j, i_3 := k$.\\

Using (\ref{omega 0 form}) and the formulae (\ref{phi -1 projection}) and (\ref{phi -1 projection a}), one sees that $$(\varphi_{-1} \circ \sigma^*\omega_{\leq})(ug_I)(\xi_J^h(ug_I)) = (\varphi_{-1} \circ \sigma^*\omega_{\leq})(ug_I)((R_{g_I})_*(\xi_J^h(u))) = (0,\tilde{\mathrm{s}}_g,0,-1)$$ and $(\varphi_{-1} \circ \sigma^*\omega_{\leq})(ug_I)(\xi_K^h(ug_I)) = (0,-\tilde{\mathrm{s}}_g i,0,i)$. Thus, if we denote by $J_0$ the complex structure on $\tilde{\g}_{-1}$ given by component-wise multiplication by $-i$, then one computes:
\begin{align*}
J_0((\varphi_{-1} \circ \sigma^*\omega_{\leq})(\xi_J^h(ug_I))) = (\varphi_{-1} \circ \sigma^*\omega_{\leq})(\xi_K^h(ug_I)); \\
J_0((\varphi_{-1} \circ \sigma^*\omega_{\leq})(\xi_K^h(ug_I))) = -(\varphi_{-1} \circ \sigma^*\omega_{\leq})(\xi_J^h(ug_I)).
\end{align*}
\no So on $\tilde{\Dbdle}_I$, the restriction of $\tilde{J}$ to $(\xi_I^{\perp})^{\nabla}$ is given on basis vectors by: $$\tilde{J}: \xi_J^{\nabla} \mapsto \xi_K^{\nabla} = (\xi_I \times \xi_J)^{\nabla} = (\xi_{I \circ J})^{\nabla} \,\,\,\, \mathrm{and} \,\,\,\, \tilde{J}: \xi_K^{\nabla} \mapsto -\xi_J^{\nabla} = (\xi_I \times \xi_K)^{\nabla} = (\xi_{I \circ K})^{\nabla},$$ where ``$\times$'' denotes the cross product in $\Vbdle_x \isom \R^3$. In a similar way, one sees that the complex structure $J_0$ on $\tilde{\g}_{-1}$ induces the natural complex structure on $Ver_I{\mathcal{Z}}$: $\tilde{J}: J \mapsto I \circ J = K$ and $\tilde{J}: K \mapsto I \circ K = -J$ ($J, K$ are naturally identified with vectors in $T_I(\mathcal{Z}_x \isom S^2)$ since they are orthogonal to $I$, and one calculates that this $\tilde{J}$ is induced from the transformation on $Ver_{ug_I}\overline{\Gbdle}_0$ which sends the fundamental vector field of $[(0,j,0)] \in \g_0$ to the fundamental vector field of $[(0,k,0)] \in \g_0$ and sending $\widetilde{[(0,k,0)]}$ to $-\widetilde{[(0,j,0)]}$).\\

Finally, to see the restriction of $\tilde{J}$ to $(\Dbdle_x)^{\nabla}$, let $X \in \Dbdle_x$ be such that $\omega_{-1}(u)(X^h) = [u]_{-1}(X) = [\overline{x}]_{-1} \in \g_{-1}$. Then $J_0(\varphi_{-1}(\omega_{-1}(u)(X^h))) = \varphi_{-1}([\overline{x}]_{-1})\overline{i} = \varphi_{-1}([\overline{x} \, \overline{i}]_{-1}) = \varphi_{-1}([\overline{ix}]_{-1})$, since $\varphi_{-1}$ is clearly $\C$-linear. On the other hand, we have $\omega_{-1}(u)(I_1(X)^h) = [u]_{-1}(I_1(X)) = [\overline{ix}]_{-1}$ (cf. Section 2.2. and Appendix A of \cite{qcweyl}), which shows that at the point $I_1 \in \mathcal{Z}_x$, the restriction of $\tilde{J}$ to $(\Dbdle_x)^{\nabla} \subset \tilde{\Dbdle}_{I_1}$ is given by $I_1^{\nabla}$. On the other hand, from the equivariance of $\omega_{-1}$ it follows that $\omega_{-1}(ug_I)(I(X)^h) = [ug_I]_{-1}(I(X)) = [ug_I]_{-1}(X)\overline{i}$. But we have $J_0(\varphi_{-1}(\omega_{-1}(ug_I)(X^h))) = \varphi_{-1}([ug_I]_{-1}(X)\overline{i})$, which shows that the restriction of $\tilde{J}$ to $(\Dbdle_x)^{\nabla} \subset \tilde{\Dbdle} \subset T_{I}\mathcal{Z}$ is also given by $I^{\nabla}$ for arbitrary $I \in \mathcal{Z}_x$. This completes the proof of Theorem A.\\

\no \begin{small}SCHOOL OF MATHEMATICS, UNIVERSITY OF THE WITWATERSRAND, P O WITS 2050, JOHANNESBURG, SOUTH AFRICA.\end{small}\\
\no E-mail: \verb"jesse.alt@wits.ac.za"

\end{document}